\documentclass[a4paper,11pt,twoside, reqno]{amsart}
\usepackage{fancyhdr,indentfirst}
\usepackage{a4wide,amssymb, amsmath, amsthm}
\usepackage[latin1]{inputenc}
\usepackage[bookmarks=false,linkcolor=blue,colorlinks=true,dvips=true]{hyperref}

\usepackage{enumitem, array}
\usepackage[round]{natbib}
\def\cite{\citet}

\emergencystretch=\hsize
\newcommand{\cc}{{\mathcal C}}

\newcommand{\cf}{{\mathcal F}}

\newcommand{\ck}{{\mathcal K}}

\newcommand{\E}{{\mathbb E}}

\renewcommand{\L}{{\mathbb L}}

\renewcommand{\P}{{\mathbb P}}

\newcommand{\R}{{\mathbb R}}

\newcommand{\s}{\star}

\DeclareMathAlphabet{\mathbbold}{U}{bbold}{m}{n}
\DeclareMathAlphabet{\mathbbmss}{U}{bbmss}{m}{n}
\SetMathAlphabet\mathbbmss{bold}{U}{bbmss}{bx}{n}
\newcommand{\ind}[1]{{\mathbbold 1}_{\{#1\}}}

\newcommand{\norm}[1]{\mathop{\left\| #1 \right\|}\nolimits}

\newcommand{\interior}[1]{\mathop {\mathrm{int}(#1)}}

\newcommand{\as}{almost surely}
\theoremstyle{plain}
\newtheorem{thm}{Theorem}

\newtheorem{lemma}{Lemma}
\newtheorem{rem}{Remark}

\theoremstyle{definition}
\newtheorem*{ack}{Acknowledgement}

\makeatletter
\newcounter{hypo}
\newcommand*{\dohypo}{\textbf{(${\mathcal A}$\thehypo)}}
\newenvironment{hypo}[1][]{%
  
  \refstepcounter{hypo}
  \list{}{%
    \settowidth{\labelwidth}{\dohypo}%
    \setlength{\labelsep}{10pt}%
    \setlength{\leftmargin}{\labelwidth}
    \advance\leftmargin\labelsep%
  }%
\item[\dohypo  #1]%
}{%
  \endlist
}

\def\hypref#1{\hyperref[hyp:#1]{(${\mathcal A}$\ref*{hyp:#1})}}
\def\hypreff#1#2{\hyperref[hyp:#2]{(${\mathcal A}$\ref*{hyp:#1}-{\it \ref*{hyp:#2})}}}

\makeatother


\author{ J\'{e}r\^{o}me LELONG}
\address{CERMICS, Ecole des Ponts, ParisTech,
    6-8 avenue Blaise Pascal, Champs sur Marne 77455 Marne La Vallée, FRANCE.}
    \email{lelong@cermics.enpc.fr}
\date{\today}
\title[A.s. convergence of randomly truncated stochastic
  algorithms]{Almost sure convergence of randomly truncated stochastic
  algorithms under verifiable conditions}

\begin{document}

\begin{abstract} In this paper, we are interested in the almost sure convergence
  of randomly truncated stochastic algorithms. In their pioneer work,
  \cite{chen86:_stoch_approx_proced}  required that the family of the noise
  terms is summable to ensure the convergence.  In our paper, we present a new
  convergence theorem which extends the already known results by making vanish
  this condition on the noise terms \textemdash\ a condition which is quite hard
  to check in practice.  The aim of this work is to prove an almost sure
  convergence result of randomly truncated stochastic algorithms under easily
  verifiable conditions (see Theorem~\ref{thm:cor}).
  \\
  
\end{abstract}

\maketitle

\section{Introduction}

The localisation of the zeros of a function $u$ is a quite complicated problem
for which many techniques have already been developed. The use of stochastic
algorithms is widely spread for solving such problems. Stochastic algorithms are
particularly well suited where some on-line parameter estimation is needed. Such
algorithms go back to the pioneer work of \cite{MR0042668}. They proposed to
consider the following recurrence relation
\begin{equation*}
  X_{n+1} = X_n - \gamma_{n+1} u(X_n) - \gamma_{n+1} \delta M_{n+1},
\end{equation*}
where $\gamma_n$ is a decreasing gain sequence and $\Delta m_n$ the measurement
error. Under certain conditions on the growth of the $\L^2-$norm of the error,
$X_n$ converges almost surely to the unique root of $u$. Since their work, much
attention has been drawn to the study of the theory of such recursive
approximations. The first works were dealing with independent measurement error
on the observations. A great effort was made in this direction to weaken the
conditions imposed on both the regressive function and the noise term.  Using
the ordinary differential equation technique, \cite{MR499560} proved a
convergence result for a wider range of measurement noises and in particular for
martingale increments.

One major drawback of these algorithms is that their convergence can only be
established if the function $u$ does not grow too quickly, namely a sub-linear
behaviour is required. This is a dramatic restriction for practical
applications.  \cite{chen86:_stoch_approx_proced} have found a way to get round
the restriction by considering stochastic algorithms truncated at randomly
varying bounds. Their algorithm can be written
\begin{equation}
  \label{eq:chen-general}
  X_{n+1} = X_n - \gamma_{n+1} u(X_n) - \gamma_{n+1} \delta M_{n+1} - \gamma_{n+1}
  p_{n+1}, 
\end{equation}
where $p_n$ is a truncation term.

In this paper, we are concerned with the convergence of the truncated
algorithm~(\ref{eq:chen-general}). Several results already exist but the
hypotheses considered differ quite significantly. The first result concerning
the almost sure convergence is due to \cite{chen86:_stoch_approx_proced}. The
convergence was also studied by  \cite{delyon96:_gener_resul_conver_stoch_algor}
and \cite{MR1701103}.  The robustness of the procedure was established by
\cite{chen88:_conver_robbin_monro} under global hypotheses on the measurement
error. Namely, they require that the series $\sum_n \gamma_n \delta M_n$
converges almost surely.  \cite{delyon96:_gener_resul_conver_stoch_algor} has
also studied the almost sure convergence under local hypotheses on the
measurement noise.  Here, we give a self-contained proof of the convergence
under local hypotheses, that is we only assume that $\sum_n \gamma_n \delta M_n
\ind{X_n \in K}$ converges almost surely for any compact set $K$. We do not
impose any condition on the truncation term $p_n$.

First, we define the general framework and explain the algorithm developed by
\cite{chen86:_stoch_approx_proced}. Our main result is stated in
Theorem~\ref{thm:chen-LF} in a very general way. For practical purposes, we give
in Theorem~\ref{thm:cor} an easily verifiable condition under which our main
result holds. This theorem is extremely valuable and dramatically extends the
range of applications of randomly truncated stochastic algorithms. Finally,
Section~\ref{sec:proof-as} is devoted to the proof of the general convergence
theorem.

\section{General framework}

Let us consider a general problem consisting in finding the root of a continuous
function $u \colon X \in \R^d \longmapsto u(X) \in \R^d$, defined as an
expectation on a probability space $(\Omega, \mathcal{A}, \P)$.
\begin{equation}
  \label{eq:def_u}
  u(X) = \E(U(X,Z)),
\end{equation}
where Z is a random variable in $\R^m$ and $U$ a measurable function defined on $\R^d
\times \R^m$ into $\R^d$. We assume that $x \longrightarrow \E(\norm{U(x,Z)}^2)$
grow faster that $\norm{x}^2$, so that the convergence of the standard Robbins Monro
algorithm is not guarantied. Instead, we consider the alternative procedure proposed
by \cite{chen86:_stoch_approx_proced}, on which we concentrate in this work.

The technique consists in forcing the algorithm to remain in an increasing sequence
of compact sets. Somehow, it prevents the algorithm from blowing up during the "first"
steps.

We consider an increasing sequence  of compact sets ${({\mathcal K}_j)}_j$
\begin{equation}
  \label{eq:chen-compact-def}
  \bigcup_{j=0}^{\infty} \mathcal{K}_j \: = \: \R^d \quad \mbox{and} \quad \forall j,
  \; \ck_j \varsubsetneq \interior{\ck_{j+1}}.
\end{equation}
We also introduce ${(Z_n)}_n$ an independent and identically distributed sequence of
random variables following the law of $Z$ and ${(\gamma_n)}_n$ a decreasing sequence
of positive real numbers. $\gamma_n$ is often called the gain sequence. For any
deterministic $X_0 \in \ck_0$ and $\sigma_0 = 0$, we define the sequences of random
variables ${(X_n)}_n$ and ${(\sigma_n)}_n$.
\begin{equation}
  \label{eq:1}
  \begin{cases}
    & X_{n + \frac{1}{2}}  = X_{n} - \gamma_{n+1} U(X_n, Z_{n+1}),\\
    \text{if $X_{n + \frac{1}{2}} \in \mathcal{K}_{\sigma_n}$} & X_{n+1} =
    X_{n + \frac{1}{2}} \quad \mbox{ and } \quad \sigma_{n+1} = \sigma_n, \\
    \text{if $X_{n + \frac{1}{2}} \notin \mathcal{K}_{\sigma_n}$} & X_{n+1}
    = X_{0}  \quad \mbox{ and } \quad \sigma_{n+1} = \sigma_n + 1. 
  \end{cases}
\end{equation}

\begin{rem} When $X_{n + \frac{1}{2}} \notin \mathcal{K}_{\sigma_n}$, one can
  set $X_{n+1}$ to any measurable function of $(X_0, \dots, X_n)$ with
  values in a given compact set. This existence of such a compact set is definitely
  essential to prove the a.s. convergence of ${(X_n)}_n$.
\end{rem}

\begin{rem} $X_{n+ \frac{1}{2}}$ represents the iterate of the Robbins Monro
  algorithm at step $n+1$.
\end{rem}

We introduce $\cf_n = \sigma(Z_k; k\leq n)$ the $\sigma$-field generated by the
random vectors $Z_k$, for $k \leq n$. Note that $X_n$ is $\cf_n-$measurable since
$X_0$ is deterministic and $U$ measurable. We can write $u(X_n) =
\E[U(X_n, Z_{n+1})|\cf_n]$.

It is often more convenient to rewrite (\ref{eq:1}) as follows
\begin{equation}
  \label{eq:theta-chen}
  X_{n+1} = X_{n}  - \gamma_{n+1} u(X_n) -
  \gamma_{n+1}\delta M_{n+1} + \gamma_{n+1} p_{n+1}  
\end{equation}
where
\begin{eqnarray}
  \label{eq:delta-Mn-2}
  \delta M_{n+1} & = & U(X_n, Z_{n+1}) -u(X_n), \\
  \label{eq:rn}
  \mbox{and} \quad  p_{n+1} & = &
  \begin{cases}
    u(X_n) + \delta M_{n+1} + \frac{1}{\gamma_{n+1}} (X_0 - X_n) &
    \text{if } X_{n+\frac{1}{2}} \notin \ck_{\sigma_n},\\ 
    0 & \text{otherwise.}
  \end{cases}
\end{eqnarray}

\begin{rem}
  $\delta M_n$ is a martingale increment. The case of the standard Robbins Monro
  algorithm corresponds to $p_n=0$.
\end{rem}

\section{Almost sure convergence}
\label{sec:as}

In this section, we present a new convergence theorem that improves the result
of \cite{chen86:_stoch_approx_proced} who proved the almost sure convergence
under global hypotheses on the series $\sum_n \gamma_{n+1} \delta M_{n+1}$
whereas we can manage the proof under local hypotheses only, namely we only
assume that the function $x \longmapsto \E(\norm{U(x,Z)}^2)$ is bounded on all
compact sets. Such a local hypothesis is much easier to satisfy in practical
applications.

\begin{thm}
  \label{thm:cor} We assume that 
  \begin{hypo}
    \label{hyp:lyapounov} There exists a unique $x^\star$ s.t. $u(x^\star)=0$
    and  $\forall \, x \neq x^\star$, $(u(x) | (x-x^\star)) > 0$.
  \end{hypo}
  \begin{hypo}
    \label{hyp:step-square}
    $\sum_n \gamma_n = \infty$ and $\sum_n \gamma_n^2 < \infty$.
  \end{hypo}
  \begin{hypo}
    \label{hyp:cv-incr-2} The function $x \longmapsto
    \E(\norm{U(x,Z)}^2)$ is bounded on any compact sets.
  \end{hypo}
  Then, the sequence ${(X_n)}_n$ converges a.s. to $x^{\star}$ for
  any sequence of compact sets satisfying~(\ref{eq:chen-compact-def}) and
  moreover the sequence ${(\sigma_n)}_n$ is a.s. finite (i.e.  for $n$ large
  enough $p_n=0$ a.s.).
\end{thm}

We will not prove Theorem~\ref{thm:cor} directly as it actually derives from a
more general result.

\begin{thm}
  \label{thm:chen-LF}
  Under Hypothesis~\hypref{lyapounov} and if
  \begin{hypo}
    \label{hyp:step}
    $\sum_n \gamma_n = \infty$.
  \end{hypo}
  \begin{hypo}
    \label{hyp:cv-incr} For all $q>0$, the series $\sum_n \gamma_{n+1} \delta M_{n+1}
    \ind{\norm{X_n - x^\s} \leq q}$ converges \as.
  \end{hypo}
  Then, the sequence ${(X_n)}_n$ converges a.s. to $x^{\star}$ and
  moreover the sequence ${(\sigma_n)}_n$ is a.s. finite (i.e.  for $n$ large
  enough $p_n=0$ a.s.).
\end{thm}

\begin{rem} In the case where $u$ derives from a potential $V$ (i.e. $u=\nabla V$),
  Hypothesis~\hypref{lyapounov} is satisfied as soon as $V$ is strictly convex.
\end{rem}

\begin{proof}[Proof of Theorem~\ref{thm:cor}] It is sufficient to prove that the
  hypotheses of Theorem~\ref{thm:cor} imply the ones of
  Theorem~\ref{thm:chen-LF}. Consider $M_n = \sum_{i=1}^{n} \gamma_{i} \delta
  M_{i} \ind{\norm{X_{i-1}-x^\s} \leq q}$, $(M_n)_n$ is a martingale. By
  computing its angle bracket, we find $\langle M \rangle_n =\sum_{i=1}^{n}
  \gamma_{i}^2 \E(\delta M_{i}\delta M_{i}^\prime | \cf_{i-1})
  \ind{\norm{X_{i-1}-x^\s} \leq q}$. As the series $\sum_i \gamma_i^2$ converges
  and the function $x \longmapsto \E(\norm{U(x,Z)}^2)$ is bounded on all compact
  sets, the almost sure convergence of $\langle M \rangle_n$ ensues from the
  Strong Law for square integrable martingales. Hence, we can apply
  Theorem~\ref{thm:chen-LF}, and the conclusion yields.
\end{proof}

\section{Proof of Theorem~\ref{thm:chen-LF}}
\label{sec:proof-as}

The proof of Theorem~\ref{thm:chen-LF} is based on the following lemma which
establishes a condition for the sequence ${(X_n)}_n$ to be a.s. compact.

\begin{lemma}
  \label{lemma:compact-sequence} If for all $q>0$, the series $\sum_{n>0} \gamma_n
  \delta M_n \ind{\norm{X_{n-1}-x^\s}<q}$ converges a.s. and if $p_n
  \ind{\norm{X_{n-1}-x^\s}<q} \longrightarrow 0$, then the sequence  
  ${(X_n)}_n$ remains a.s. in a compact set.
\end{lemma}

Note that the compact set mentioned in Lemma~\ref{lemma:compact-sequence} is
random.  In particular, this lemma does not imply that the number of truncations
is bounded independently of the randomness~$\omega$.

\begin{proof}[Proof of Theorem~\ref{thm:chen-LF}]
  The proof is divided in two parts. 
  \begin{itemize}
  \item  Let $q>0$. We define $\bar M_n = \sum_{i=1}^n \gamma_i \delta M_i
    \ind{\norm{X_{i-1}-x^\s} \leq q }$. Thanks to Hypothesis~\hypref{cv-incr}, $\bar
    M_n$ converges \as.

    Assume that $\sigma_n \longrightarrow \infty$. This is in contradiction with
    the conclusion of Lemma~\ref{lemma:compact-sequence}, which implies that the
    hypothesis according to which $p_n \ind{\norm{X_{n-1}-x^\s}<q}$ tends to $0$
    does not hold. So,
    \begin{equation*}
      \exists \: \eta>0,\: q>0,  \quad \forall N>0, \: \exists n>N \quad
      \ind{\norm{X_n -x^\s} \leq q} \norm{p_{n+1}} > \eta.
    \end{equation*}

    Let $\varepsilon>0$. There exists a subsequence $X_{\phi(n)}$ such that for all
    $n>0$, $\ind{\norm{X_{\phi(n)}-x^\s} \leq q} \norm{p_{{\phi(n)}+1}} \neq 0$ and
    $\norm{\gamma_{{\phi(n)}+1} \delta M_{{\phi(n)}+1}} \leq \varepsilon$.

    So, $\norm{X_{\phi(n)}-x^\s} \leq q$ and however the new potential iterate
    $X_{{\phi(n)}+\frac{1}{2}} = X_{\phi(n)}  - \gamma_{{\phi(n)}+1} (u(X_{\phi(n)})
    + \delta M_{{\phi(n)}+1})$ is not in $\ck_{\sigma_{\phi(n)}}$. Since $u$ is
    continuous, $\norm{\gamma_{{\phi(n)}+1} u(X_{\phi(n)})}$ can be made smaller than
    $\varepsilon$. As $\norm{\gamma_{{\phi(n)}+1} \delta M_{{\phi(n)}+1}} \leq
    \varepsilon$, a proper choice of $\varepsilon$ enables to write
    \begin{equation*}
      \norm{X_{\phi(n)} - x^\s  - \gamma_{{\phi(n)}+1}
        (u(X_{\phi(n)}) + \delta M_{{\phi(n)}+1})} \leq q+1.
    \end{equation*}
    Let $l$ be the smallest integer s.t. $B(x^\s, q+1) \subset \ck_l$ (such an integer
    exists thanks to (\ref{eq:chen-compact-def})), then $\sigma_{\phi(n)} < l$ for
    all $n$. Since the sequence $(\sigma_n)_n$ is increasing, this proves that
    $\limsup_n \sigma_n < \infty$ a.s..

  \item According to the previous item $\limsup_n \sigma_n < \infty$ a.s.. So, the
    sequence $(X_n)_n$ is \as\ compact. Consequently, we can in fact set $q = \infty$
    in Hypothesis~\hypref{cv-incr} and say that $\sum_{i} \gamma_i \delta M_i$
    converges \as. Let us consider
    \begin{equation*}
      X_n^\prime = X_n - \sum_{i=n+1}^\infty \gamma_i \delta M_i.
    \end{equation*}

    Since the series $\sum_{i>0} \gamma_i \delta M_i$ converges a.s. and $X_n$ remains
    in a compact set, $X_n^\prime$ also remains in a  compact set. Let $\cc$ be this
    compact set. We define $\bar u = \sup_{x \in \cc} \norm{u(x)}$.

    \begin{equation*}
      X_{n+1}^\prime = X_n^\prime - \gamma_{n+1} u(X_n^\prime) +
      \gamma_{n+1} \varepsilon_n,
    \end{equation*}
    where $\varepsilon_n = u(X_n^\prime) - u(X_n)$. Since $\norm{X^\prime_n -X_n}
    \longrightarrow 0$ and $u$ is continuous, $\norm{\varepsilon_n} \longrightarrow
    0$.  
    \begin{multline*}
      \qquad \norm{X_{n+1}^\prime - x^\star}^2  \leq  \norm{X_n^\prime -
        x^\star}^2 - 2 \gamma_{n+1} (X_n^\prime -
      x^\star \: | \: u(X_n^\prime)) \\
      + \gamma_{n+1}^2 (\varepsilon_n^2 +
      \bar u^2) -2  \gamma_{n+1} (X_n^\prime -
      x^\star \: | \: \varepsilon_n).
    \end{multline*}
    We can rewrite the inequality introducing a new sequence
    $\varepsilon_n^\prime \longrightarrow 0$.
    \begin{eqnarray}
      \label{eq:30}
      \norm{X_{n+1}^\prime - x^\star}^2  & \leq  & \norm{X_n^\prime -
        x^\star}^2 - 2 \gamma_{n+1} (X_n^\prime -
      x^\star \: | \: u(X_n^\prime)) + \gamma_{n+1} \varepsilon_n^\prime.
    \end{eqnarray}

    Let $\delta>0$. If $\norm{X_n^\prime - x^\star}^2 > \delta$, then $(X_n^\prime -
    x^\star \: , \: u(X_n^\prime)) > c > 0$. Henceforth, for $n$ large enough
    Equation~(\ref{eq:30}) becomes
    \begin{eqnarray*}
      \norm{X_{n+1}^\prime - x^\star}^2  & \leq  & \norm{X_n^\prime -
        x^\star}^2 - \gamma_{n+1} c \ind{\norm{X_n^\prime - x^\star}^2 > \delta} +
      \gamma_{n+1} (\bar c+\varepsilon_n^\prime) \ind{\norm{X_n^\prime - x^\star}^2
        \leq \delta},
    \end{eqnarray*}
    where $\bar c = \sup_{\norm{x - x^\star}^2 \leq \delta} (x-x^\star | u(x))$.
    Since $\sum_n \gamma_n = \infty$, each time $\norm{X_n^\prime - x^\star}^2
    > \delta$, the sequence $X_n^\prime$ is driven back into the ball $\bar B(x^\star,
    \sqrt{\delta})$ in a finite number of steps. Hence, for any $n$ large enough
    \begin{equation*}
      \norm{X_n^\prime - x^\star}^2 < \delta + \gamma_{\phi(n)+1}
      (\bar c+\varepsilon_{\phi(n)}^\prime),
    \end{equation*}
    where $\phi(n) = \sup\{ p\leq n ; \norm{X'_p - x^\star}^2 \leq \delta\}$.
    As $\phi(n)$ a.s. tends to infinity with $n$, $\limsup_n \norm{X_n^\prime -
      x^\star}^2 \leq \delta$ for all $\delta>0$.  This proves that $X_n^\prime
    \longrightarrow x^\star$. Finally, since the series $\sum_n \gamma_{n+1}
    \delta M_{n+1}$ converges, this also proves that $X_n \longrightarrow
    x^\star$.
  \end{itemize}
\end{proof}

Now, we are going to prove Lemma~\ref{lemma:compact-sequence}.

\begin{proof}[Proof of Lemma~\ref{lemma:compact-sequence}] If $\sigma_n < \infty$
  a.s., the conclusion of the Lemma is obvious. Assume that  $\sigma_n
  \longrightarrow \infty$. Since each time $\sigma_n$ increases, the  sequence $X_n$
  is reset to a fixed point of $\ck_0$, the existence of a compact set in which the
  sequence lies infinitely often is straightforward.

  Let $M>0$, we set $\cc = \{x\: ; \norm{x - x^\star}^2 \leq M\}$.  We can rewrite
  the Hypotheses of the Lemma as follows
  \begin{equation}
    \label{eq:hypo-lem}
    \forall \varepsilon>0, \: \exists \: N>0 \mbox{ s.t. } \forall \:n,\:p \geq
    \:N \mbox{ we have } \; \left\{
      \begin{array}[c]{l}
        \displaystyle \norm{\sum_{k=n}^p \gamma_k \delta M_k
          \ind{ \norm{X_{k-1} - x^\star}^2 \leq M+2}} \: < \: \varepsilon, \\ 
        \gamma_n < \varepsilon,\\
        \ind{\norm{X_{n-1}-x^\star}^2 \leq M+2} \norm{p_n} < \varepsilon.
      \end{array}
    \right.
  \end{equation}
  Let $\varepsilon>0$ and $N>0$ satisfying Condition~(\ref{eq:hypo-lem})
  and s.t. $X_N \in \cc$. We introduce
  \begin{equation*}
    X_n^\prime = X_n - \sum_{i=n+1}^\infty \gamma_i \delta M_i \ind{
      \norm{X_{i-1} - x^\star}^2 \leq M+2}.
  \end{equation*}
  By using Equation~(\ref{eq:theta-chen}), we can easily show that $X_n^\prime$
  satisfies the following recurrence relation
  \begin{equation}
    \label{eq:28}
    X_{n+1}^\prime = X_n^\prime - \gamma_{n+1} \delta M_{n+1}
    \ind{\norm{X_{n} - x^\star}^2 > M+2} - \gamma_{n+1} (u(X_n) - p_{n+1}). 
  \end{equation}

  We will now prove that the sequence ${(X_n^\prime)}_n$ remains in the set $\{x \: ;
  \norm{x - x^\star}^2 \leq M+1\} \: = \: \cc^\prime$.
  
  The recurrent hypothesis is satisfied for $n = N$ (it is sufficient to choose
  $\varepsilon < 1$). Assume that the hypothesis holds for 
  $N, \dots, n$. Hence, $\norm{X_{n} - x^\star}^2 \leq M+2$. Then, we can deduce from
  Equation~(\ref{eq:28}) that
  \begin{eqnarray*}
    X_{n+1}^\prime & = & X_{n}^\prime  - \gamma_{n+1} (u(X_n) -
    p_{n+1}),\\
    \norm{X_{n+1}^\prime - x^\star}^2  & \leq  & \norm{X_n^\prime -
      x^\star}^2 - 2 \gamma_{n+1} (X_n^\prime -
    x^\star \: | \: u(X_n)) + \gamma_{n+1} c \: \varepsilon,
  \end{eqnarray*}
  where $c$ is a positive constant independent of $M$.

  \begin{itemize}
  \item   If $\norm{X_{n}^\prime - x^\star}^2 \leq M$, thanks to the continuity of
    $u$, a proper choice of $\varepsilon$ ensures that $\gamma_{n+1} \norm{(X_n^\prime -
      x^\star \: | \: u(X_n))} < 1$. Hence, $\norm{X_{n+1}^\prime - x^\star}^2
    \leq M+1$.
  \item   If $M < \norm{X_{n}^\prime - x^\star}^2 \leq M+1$, thanks to the
    continuity of $u$ and thanks to Hypothesis~\hypref{lyapounov}, $(X_n - x^\star \: | \:
    u(X_n)) > \delta > 0$. Once again, properly choosing $\varepsilon$ guaranties that $c
    \varepsilon < \delta$. Consequently, $\norm{X_{n+1}^\prime - x^\star}^2 \leq M+1$.
  \end{itemize}

  We have proved that for all $n>N$, $\norm{X_n^\prime - x^\star}^2 \leq M+1$. Since
  $\varepsilon$ can be chosen smaller than $1$, the following upper-bound also holds
  \begin{equation*}
    \norm{X_n- x^\star}^2 \leq M+2, \mbox{ for all } n>N.    
  \end{equation*}

  This achieves to prove that the sequence ${(X_n)}_n$ remains in a compact set and
  consequently that $\limsup_n \sigma_n$ is a.s. finite.
\end{proof}

\begin{ack}
  I would like to thank Bernard Lapeyre for the fruitful remarks he
  made on a previous version of the proof presented above.
\end{ack}

\bibliographystyle{abbrvnat}
\bibliography{biblio_as}
\end{document}